\documentclass[11pt]{article}
\usepackage{amsmath,amsfonts,amssymb}
\usepackage{url}
\usepackage{hyperref}
\usepackage{srcltx}
\usepackage{xcolor}
\usepackage{epsfig}
\def\e{\mathrm{e}}
\def\mq{\mathfrak{q}}
\def\SG{{\cal SG}}
\def\oX{\overline{X}}
\def\bT{{\mathbf{T}}}
\def\wt#1{\widetilde{#1}}
\def\wtX{\widetilde{X}}

\def\mR{{{\mathfrak{R}}}}

\def\s{{\hbox{\rm\scriptsize s}}}

\def\three?{3}
\def\four?{4}
\def\ten?{10}
\def\bM{{\mathbf{M}}}

\def\bL{{\mathbf{L}}}

\def\Argmin{\mathop{\hbox{\rm Argmin}}}
\def\beq{\begin{equation}}
\def\eeq{\end{equation}}
\newtheorem{observation}{Observation}[section]

\def\norm2to2{{\|\cdot\|_{2,2}}}
\def\Prob{\hbox{\rm Prob}}

\def\cl{{\hbox{\rm  cl}\,}}

\def\bE{{\mathbf{E}}}

\def\Diag{\hbox{\rm  Diag}}
\def\Prob{\hbox{\rm  Prob}}

\def\Opt{\hbox{\rm Opt}}

\def\Conv{\hbox{\rm  Conv}}

\def\Tr{{\mathop{\hbox{\rm  Tr}}}}

\def\cD{{\cal D}}

\def\cH{{\cal H}}

\def\cN{{\cal N}}

\def\cR{{\cal R}}
\def\cS{{\cal S}}
\def\cT{{\cal T}}

\def\cW{{\cal W}}
\def\cX{{\wtX}}

\def\P{{\cal P}}

\def\U{{\cal U}}

\def\Argmin{\mathop{\hbox{\rm  Argmin}}}

\def\Ker{{\hbox{\rm  Ker}}}

\def\bL{{\mathbf{L}}}
\def\bS{{\mathbf{S}}}

\def\bH{{\mathbf{H}}}
\def\e{{\hbox{\rm e}}}

\def\qed{\ \hfill$\square$\par\smallskip}

\newtheorem{lemma}{Lemma}
\newtheorem{proposition}{Proposition}
\newtheorem{corollary}{Corollary}
\newtheorem{theorem}{Theorem}
\newtheorem{example}{Example}

\def\bR{{\mathbf{R}}}

\def\cH{{\cal H}}
\def\Col{{\hbox{\rm Col}}}

\def\Risk{{\hbox{\rm Risk}}}
\def\mP{{\mathfrak{P}}}
\def\mP{{\cal R}}
\newcommand{\be}{\begin{eqnarray}}
\newcommand{\ee}[1]{\label{#1}\end{eqnarray}}
\newcommand{\nn}{\nonumber \\}
\newcommand{\ese}{\end{align*}}
\newcommand{\bse}{\begin{align*}}
\newcommand{\rf}[1]{~(\ref{#1})}
\newcommand{\wh}[1]{{\widehat{#1}}}
\def\mR{{\mathfrak{R}}}

\newcommand{\hide}[1]{{}}

\newcommand{\anc}[2]{{\color{violet} #2}}
\title{On Design of Polyhedral Estimates in Linear Inverse Problems}
\author{
Anatoli Juditsky
\thanks{LJK, Universit\'e Grenoble Alpes, 700 Avenue Centrale,  38401 Domaine Universitaire de Saint-Martin-d'Hères, France {\tt anatoli.juditsky@univ-grenoble-alpes.fr}}
\and Arkadi Nemirovski
\thanks{Georgia Institute
 of Technology, Atlanta, Georgia
30332, USA, {\tt nemirovs@isye.gatech.edu}\newline
This work was supported by Multidisciplinary Institute in Artificial intelligence MIAI {@} Grenoble Alpes (ANR-19-P3IA-0003).}}
\date{}
\begin{document}
\maketitle
\begin{abstract}
{\em Polyhedral estimate} \cite{juditsky2020polyhedral,PUP} is a generic efficiently computable nonlinear in observations routine for recovering unknown signal belonging to a given convex compact set from noisy observation of signal's linear image. Risk analysis and optimal design of polyhedral estimates may be addressed through efficient bounding of optimal values of optimization problems. Such problems are typically hard; yet, it was shown in
\cite{juditsky2020polyhedral} that nearly minimax optimal (``up to logarithmic factors'') estimates can be efficiently constructed when the signal set is {\em an ellitope}---a member of a wide family of convex and compact sets of special geometry  (see, e.g., \cite{JudNem2018}). The subject of this paper is a new risk analysis for polyhedral estimate in the situation where the signal set is an intersection of an ellitope and an arbitrary polytope allowing for improved  polyhedral estimate design in this situation.
\end{abstract}
\section{Introduction}
In this paper we consider the estimation problem as follows. Given a noisy observation
\be
\omega&=Ax+\xi_x\in\bR^m,
\ee{eq:mod1}
of the linear image $Ax$ of unknown signal $x\in\bR^n$, known to belong to a given signal set---a nonempty convex compact set $X\subset\bR^n$, we want to recover the image $w=Bx\in\bR^\nu$ of this signal. Here $A$ and $B$ are given $m\times n$ and $\nu\times n$ matrices, and $\xi_x$ is the observation noise with distribution $P_x$ which may depend on $x$; the estimation error is measured by the $\ell_\theta$-loss  with given $\theta\in[1,2]$.

Estimation problem \rf{eq:mod1} is a  linear inverse problem.  When statistically analysed, popular approaches to solving \rf{eq:mod1} (cf., e.g., \cite{natterer1986mathematics,johnstone1990speed,johnstone1991discretization,mair1996statistical,donoho1995nonlinear,abramovich1998wavelet,goldenshluger1999pointwise,goldenshluger2000adaptive,vogel2002computational,kaipio2006statistical,cohen2004adaptive,bissantz2007convergence,hoffmann2008nonlinear,proksch2018multiscale}) usually assume a special structure of the problem, when matrix $A$ and set $X$ ``fit each other,'' e.g., there exists a sparse approximation of the set $\cX$ in a given basis/pair of bases, in which matrix $A$ is ``almost diagonal'' (see, e.g. \cite{donoho1995nonlinear,cohen2004adaptive} for detail). Under these assumptions, traditional results focus on estimation algorithms which are both numerically straightforward and statistically (asymptotically) optimal with closed form analytical description of estimates and corresponding risks.
In this paper, $A$ and $B$ are ``general'' matrices of appropriate dimensions, and $X$ is a rather general convex and compact set. Instead of deriving closed form expressions for estimates and risks, we adopt an ``operational'' approach initiated in \cite{Don95} and further developed in \cite{JN2009,l2estimation,juditsky2020polyhedral,PUP}, within which
 both the estimate and its risk are yielded by efficient computation, rather than by an explicit analytical description.

In particular, the {\em polyhedral estimate} goes back to \cite{oldpaper} where it was shown (see also \cite[Chapter 2]{saintflour}) that a polyhedral estimate is near-optimal when recovering smooth multivariate regression function known to belong to Sobolev balls from noisy observations taken along a regular grid. Recently, the ideas underlying the results of \cite{oldpaper} have been taken up in the MIND estimator of \cite{grasmair2018variational} and applied in the indirect observation setting in the context of multiple testing  in \cite{proksch2018multiscale}. The idea of this estimate, as it was reintroduced in the present setting in \cite{juditsky2020polyhedral}, may be explained as follows. Assuming that the observation noise is $\cN(0,\sigma^2I_m)$, one may evaluate linear forms $g_h^Tx=h^TAx$ of signal $x$ underlying observation \rf{eq:mod1}. When $h\in\cH=\{h\in\bR^m:\|h\|_2=1\}$, the ``plug-in'' estimate $\widehat{g}_h(\omega)=h^T\omega$ is unbiased with $\cN(0,\sigma^2)$ estimation error. When selecting somehow the matrix $H=[h_1,...,h_M]$ with columns from $\cH$, vector $H^T\omega$ approximates $H^TAx$ in the uniform norm in the sense that for any $\epsilon\in(0,1)$,
\begin{equation}\label{1Eq1}
\Prob_{\omega\sim\cN(Ax,\sigma^2I_m)}\Big\{\|H^T\omega-H^TAx\|_\infty \leq \underbrace{\sigma q_{1-\epsilon/(2M)}}_{\varrho_\epsilon}\Big\}\geq 1-\epsilon
\end{equation}
where $q_p$ is the $p$-quantile of the standard normal distribution. This estimate is then combined with a priori information that $x\in X$ to obtain the {\em polyhedral estimate}  $\widehat{w}^H(\cdot)$ of $w=Bx$ according to
\begin{equation}\label{2Eq1}
\wh{x}^H(\omega)\in\Argmin_y\left\{\|H^T[Ay-\omega]\|_\infty:\,y\in X\right\},\;\;\widehat{w}^H(\omega)=B\wh{x}^H(\omega).
\end{equation}
Similarly to  ``sketching estimators'' (cf. \cite{gribonval2021compressive} and references therein), polyhedral estimate aims to reduce
the original estimation problem to that of estimating a ``small'' set of linear forms.

Clearly, the quality of the polyhedral estimates heavily depends on the {\em design parameter}---the {\em contrast matrix} $H$. The goal of the design of a polyhedral estimate is, given matrices $A$, $B$, the set  $X$ of admissible signals and confidence parameter $\epsilon\in(0,1)$, to specify $H$ resulting in as small as possible upper bound on the $\epsilon$-risk of estimation---the minimal radius of the confidence ball for $w=Bx,\,x\in X$ of reliability $1-\epsilon$ centered at $\widehat{w}^H(\omega)$, e.g.,
$$
\Risk_{\epsilon}[\widehat{w}^H|X]:=\sup_{x\in X}\inf\left\{\rho:\Prob_{\xi_x\sim \cN(0,\sigma^2)}\{\|Bx-\widehat{w}^H(Ax+\xi_x)\|_2>\rho\}\leq \epsilon\right\}.
$$
 Note that \rf{1Eq1} implies that $\|H^TA(\widehat{w}^H(\omega)-x)\|_\infty\leq 2\rho_\epsilon$ with probability $1-\epsilon$, and because $\widehat{w}^H(\omega)\in X$ by construction, one also has $\widehat{w}^H(\omega)-x\in 2X_\s$ where $X_\s=\tfrac{1}{2} (X-X)$ is the symmetrization of $X$.
As a result, the $\epsilon$-risk of the estimate $\widehat{w}^H(\omega)$ may be bounded with the quantity (cf. \cite[Proposition 5.1]{PUP})
\[
\mR[H]=\max_{z} \big\{2\|Bz\|_2:\;z\in X_\s,\;\|H^TAz\|_\infty\leq \rho_\epsilon\big\},
\]
while the design of $H$ reduces to the optimization problem
\be\min_{H=[h_1,...,h_M],\;h_i\in \cH} \mR[H].
\ee{opt0}
Problem \rf{opt0} is generally difficult because  function $\mR[H]$ is usually nonconvex and hard to compute. However, one can replace it
with an approximating problem of minimizing w.r.t. $H$ an efficiently computable convex upper bound on $\mR[H]$. Naturally, the quality of the resulting design depends on how good the approximation in question is what, it its turn, depends on the geometric structure of $X$.
Two approaches to the design of the polyhedral estimate were proposed in \cite{juditsky2020polyhedral}.
In the first approach,  in the situation where the signal set $X$ is an {\em ellitope},\footnote{See Section \ref{sitgoal} for the formal definition; for the time being, it suffices to mention that an ellitope is a special convex and compact set which is symmetric w.r.t. to the origin, an instructive example being finite intersection of set $\|A_ix\|_{p_i}\leq1$ with $p_i\geq2$, $i=1,...,I$. Ellitopes form a rich family of convex sets which allow for tight upper-bounding of the maxima of quadratic forms using semidefinite relaxation. } the semidefinite relaxation was used to build an efficiently computable upper bound for $\mR[H]$. Subsequently, this bound  was  minimized w.r.t. to the contrast $H$. Furthermore, it was shown
\cite[Proposition 5.10]{PUP}  that in this case, the polyhedral estimate $\wh w^{H_*}$ utilizing the optimized contrast $H_*$ is nearly minimax optimal within the logarithmic
factor on the class of all estimates $\wh w$ of $w$.
The second approach implemented in \cite{juditsky2020polyhedral} was based on the observation that in the case where
the image $\cW$ of $X_\s$ under the mapping $x\mapsto Bx$ is contained in a``scaled ball'' of norm $\|\cdot\|_p$, $p\geq 1$, the quantity $\cR[H]$ admits an efficiently
computable upper bound, cf. \cite[Proposition 5.2]{PUP}. Minimizing the corresponding upper bound w.r.t. $H$ leads to the polyhedral estimate $\wh w^H$ which  is nearly minimax optimal in the
{\em diagonal case} allowing for analytical analysis \cite{donoho1994minimax,grasmair2018variational} where matrices $A$ and $B$ are diagonal and  $X=\{x\in \bR^n:\,\|Dx\|_q\leq 1\}$ with diagonal matrix $D$. Although the risk bounds for the polyhedral estimate can be efficiently upper-bounded in the ``general situation'' (when $A$, $B$ and $D$ cannot be diagonalized in the same basis), no near-optimality result is available for the corresponding polyhedral estimate in this case, and the problem of design of better contrasts in this situation is open.

An interesting feature of the polyhedral estimate \rf{2Eq1} is that it allows for a straightforward contrast aggregation. Indeed, given contrast matrices $H_k\in m\times M_k$, $k=1,...,K$, which are designed by appropriate routines, one can easily aggregate them into a single contrast matrix $H=[H_1,...,H_K]$, resulting in the polyhedral estimate $\wh w^H$ with the maximal risk which is similar to $\min_k \cR_k$ up to a moderate (logarithmic in $K$) factor. This property of the estimate allows, for instance, for aggregating contrasts obtained when implementing the design approaches discussed above, so that the risk bound of the ``aggregated'' estimate is, within log-factor,  the minimum of the risk bounds obtained by the corresponding approximations.
For instance, when the signal set $X$ is an intersection of two sets, say, an ellitope $X'$ and a ball of $\ell_1$-norm $X''$, this property of the estimate allows to aggregate the contrasts $H'$ and $H''$ designed for each signal set when utilizing the two approaches  discussed above, with the risk bound of the ``synthetic'' estimate being nearly the minimum of the respective risk bounds derived in the two approaches. Note, that when the signal set $X'\cap X''$ is ``significantly smaller'' than each of the sets, one would expect the corresponding risk bound to be better than the simple minimum of two risks. However, this cannot be achieved when implementing the aggregation routine in question.

In this work, our goal is to overcome to some extent the mentioned above drawbacks of the polyhedral estimate design as presented in \cite{juditsky2020polyhedral} and \cite[Section 5.1]{PUP}.
Specifically, we propose new bounds for the quantity $\cR[H]$ when the symmetrization of the signal set $X$ is the intersection of a given ellitope and the convex hull of finitely many given vectors. Although based on the same ideas as the second bounding scheme of \cite{juditsky2020polyhedral} mentioned above, new bounds lead to the contrast design resulting in the improved accuracy of the polyhedral estimate in this situation.
As a part of our developments, we introduce a simple but seemingly new risk decomposition procedure for more precise risk bounding in the problem of estimation over intersection of the sets of different geometry, specifically, a polytope and an ellitope.

The paper is organized as follows. In Section \ref{sec:2} we describe the estimation problem setting and recall some properties of the polyhedral estimate  to be used in the sequel. We describe our new approach to bounding the estimation risk and contrast design for building a ``presumably good'' polyhedral estimate in Section \ref{secdesign}. Finally, we report on some preliminary numerical experiments in Section \ref{sec:num}.

\par
Proofs of the results are postponed till the appendix.
\section{Situation and goal}\label{sec:2}
\subsection{Problem statement}\label{sitgoal}
Recall that, by definition \cite{JudNem2018,PUP}, {\em an ellitope} in $\bR^k$ is a set of the form
\[   Y=\{y\in\bR^k:\,\exists z\in\bR^K,\,s\in\cS:\,y=Sz,\,z^TS_\ell z\leq s_\ell,\,\ell\leq L\},
\] where $S_\ell\succeq0$, $\sum_\ell S_\ell \succ0$, and $\cS\subset\bR^K_+$ is a convex compact set with a nonempty interior which is monotone: whenever $0\leq s'\leq s\in\cS$ one has $s'\in\cS$.
We refer to $L$ as {\em ellitopic dimension} of $Y$.\footnote{
Ellitopes form a large family of symmetric w.r.t. the origin convex and compact sets. This family is closed w.r.t. some basic operations preserving convexity, compactness, and symmetry such as taking finite intersections, direct products, linear images, etc., and  also
allows for an algorithmic ``calculus'' which applies to ``raw materials,'' e.g., $\|\cdot\|_p$-balls with $p\in[2,\infty]$.}

From now on we assume that the convex and compact signal  set $X\subset \bR^n$ is such that its symmetrization  $X_\s={1\over 2}[X-X]$ belongs to the intersection of an ellitope
\beq\label{cX}
\wt{X}=\left\{x\in\bR^n:\;\exists z\in \bR^N,\, t\in \cT:\, x=Pz, z^TT_kz\leq t_k,k\leq K\right\}
\eeq
with the set
\beq\label{ovX}
\oX=\left\{x\in\bR^n:\,\exists w\in\bR^q:\,x=Rw,\,Qw\in\Conv\{\pm v_1,...,\pm v_J\}\right\}
\eeq
where $p\times q$ matrix $Q$ is with trivial kernel, $R$ is $n\times q$ matrix, and $v_j\in\bR^p$.

We consider signal recovery problem as follows: given an observation (cf. \rf{eq:mod1})
$$
\omega=Ax+\xi_x\in \bR^m,\quad A\in\bR^{m\times n}, \quad \xi_x\sim P_x,
$$
our objective is to estimate the linear image $w=Bx$, $B\in \bR^{\nu\times n}$, of unknown signal $x$ known to belong to  $X$.

Given $\epsilon\in(0,1)$ and once for ever fixed $\theta\in[1,2]$, we quantify the quality of  recovery $\wh w(\cdot)$ by its {\em $\epsilon$-risk}
\[
\Risk_{\epsilon,\theta}[\wh{w}|X]:=\sup_{x\in X}\inf\left\{\rho:\,\Prob_{\xi\sim P_x}\{\|w-\widehat{w}(\omega)\|_\theta>\rho\}\leq \epsilon\right\}.
\]
We assume that the distribution $P_x$ of the observation noise $\xi_x$ is fully specified by $x$. Furthermore, for every $\delta>0$ we have at our disposal a norm $\pi_\delta(\cdot)$ such that whenever $\pi_\delta(h)\leq1$, one has
\[
\Prob_{\xi\sim P_x}\{|h^T\xi|>1\}\leq \delta\;\forall x\in X.
\]
Besides this, we suppose that the unit ball of the norm $\pi_\delta(\cdot)$ is an ellitope:
\begin{align}
B^{\delta}&:=\{h\in\bR^m:\,\pi_\delta(h)\leq1\}\nn
&=\left\{h:\,\exists z\in\bR^{M},\,s\in\cS^\delta:\, h=S_\delta z,\,z^TS^\delta_\ell z\leq s_\ell,\,\ell\leq L\right\}.
\label{Bpi}
\end{align}
We suppose that $L$ and $M$, in contrast to other components in the ellitopic description of $B^{\delta}$, are independent of $\delta$.
Note that the simplest example of this situation is the zero mean Gaussian observation noise s$\xi_x\sim\cN(0,\sigma_x^2I_m)$, $\sigma_x\leq\overline{\sigma}$, $x\in X$; in this case we can take
$$
B^\delta=\{h\in\bR^m:\exists s\in [0,\overline{\sigma}^{-2}q^{-2}_{1-\delta/2}]:h^Th\leq s\},
$$
where $q_\gamma$ is $\gamma$-quantile of $\cN(0,1)$, that is, $\pi_\delta(h)=\overline{\sigma}q_{1-\delta/2}\|h\|_2$.
\par
Finally, we say that matrix $H$ of column dimension $m$ is {\em $\delta$-admissible} if $\pi_\delta$-norms of the columns of $H$ do not exceed 1.

\subsection{Application: signal recovery in mixture-sub-Gaussian regression}\label{sec:numeric}
To motivate the proposed problem setting, consider a particular application in which the parameter to be recovered is naturally localized in a polytope.  Specifically, consider the following situation:
there is a stream of particles of $n$ different types arriving at a detecting device. At time instant $t=1,2,...$ the detector is hit by a particle of type $i\in \{1,...,n\}$, with $i$  drawn at random according to the proportions $x_j$ of particles of different types in the stream. The detector's output $\omega_t\in\bR^d$ is the ``signature'' of the particle drawn at random from some distribution $P_t$ from a given family of distributions, and our goal is to infer from $N$ observations the distribution  $x$ of particles of different types.
\par
Consider the following model of the above situation: suppose that we are given $n$ pairs $a_i\in\bR^d$, $\Theta_i\in\bS^d_+$ which define $n$ families $\SG(a_i,\Theta_i)$ of sub-Gaussian
probability distributions on $\bR^d$ with parameters $a_i,\;\Theta_i$.\footnote{We say that random vector $\eta\in \bR^k$ is sub-Gaussian with parameters $\mu$ and $\Sigma$ ($\eta\sim \SG(\mu, \Sigma))$ if $\bE\big\{e^{u^T\eta}\big\}\leq \exp\Big\{u^T\mu+\tfrac{1}{2} u^T\Sigma u\Big\}$  $\forall u\in \bR^k$.}
At time $t$, $1\leq t\leq N$, a realization $\omega_t$ is generated as follows:  nature
draws at random index $i_t\in I=\{1,..,n\}$ according to the probability distribution $x$ on $I$, and then selects a probability distribution $P_{i_t}$ from $\SG(a_{i_t},\Theta_{i_t})$
and draws $\omega_t$ from this distribution. The available observation is the mean of $\omega_t$ over $t=1,...,N$,
\[
\omega={1\over N}\sum_{t=1}^N \omega_t,
\]
and our goal is to recover $x$.
We assume that nature actions are independent across time, so that $\omega_t$, $t=1,...,N$, are independent but not necessarily identically distributed.
\par
To comply with the framework of Section \ref{sitgoal}, we specify $A$ as the $d\times n$ matrix with columns $a_1,...,a_n$. Given a vector  $x\in\bR^n_+$ with unit sum of entries, let us set
$$
a=Ax,\;b_i=a_i-a,i\leq n,\;\xi^t_x=\omega_t-a,\;\xi_x={1\over N}\sum_{t=1}^N\xi^t_x,
$$
thus arriving at the observation scheme
$$
\omega=Ax +\xi_x,
$$
as required by the setup of Section \ref{sitgoal}. We assume that $X$ is the intersection of the standard simplex $\Big\{x\in\bR^n_+:\,\sum_ix_i=1\Big\}$ with an ellitope $\cX$, so that
\[
X_\s=\cX\cap \oX,\;\;\oX=\Conv\{\pm v_1,\pm v_2,...,\pm v_n\},
\]
where $v_j$ are the standard basis vectors in $\bR^n$. To complete the problem specifications, we need to point out the norm $\pi_\delta(\cdot)$ we intend to use in our construction.
We utilize the following simple
fact:
\begin{lemma}\label{lem1} Let $h\in\bR^d$ satisfy
\[
|h^T[a_i-a_j]|\leq\varkappa,\;h^T\Theta_ih\leq \sigma^2,\;\;1\leq i,j\leq n.
\]
Then for every $x\in X$ it holds
\be
 \Prob\big\{|h^T\xi_x|>\tau \big\}\leq 2\exp\left\{-\frac{\tau^2N}{2(\varkappa^2+\sigma^2)}\right\},\;\; \forall \tau\geq 0.
\ee{prob000}
 Given $\delta\in(0,1)$ and a positive integer $N$, consider the norm
\[
\pi_\delta(g)=2\beta^{-1}\max\left[\max_{i,j}|g^T[a_i-a_j]|,\;\max_i\sqrt{g^T\Theta_ig}\right],\;\;\beta=\sqrt{N\over \ln(2/\delta)}.
\]
Then for every $x\in X$ and every $g\in\bR^n$ one has
\beq\label{ifthen}
\pi_\delta(g)\leq 1\Rightarrow \Prob\{|g^T\xi_x|>1\}\leq\delta.
\eeq
\end{lemma}
By Lemma \ref{lem1}, to fit all assumptions in the setup of Section \ref{sitgoal} it suffices to specify
the missing so far setups's components---the entities  participating in (\ref{Bpi})---according to
\begin{align*}
L&=n(n+1)/2,\;
M=n,\\
\cS^\delta&=\{s=\{s_{ij},1\leq i\leq j\leq n\}:0\leq s_{ij}\leq1,0\leq i\leq j\leq n\} \\
S^\delta_{ij}&=
\tfrac{4}{\beta^2}[a_i-a_j][a_i-a_j]^T,\,1\leq i<j\leq n,\;\;
S^\delta_{ii}=\tfrac{4}{\beta^2}\Theta_i,1\leq i\leq n.
\end{align*}

\subsection{Polyhedral estimate: preliminaries}
Under the circumstances, a {\sl polyhedral estimate} is specified by $\delta\in(0,1)$ and $m\times \mu$  $\delta$-admissible contrast matrix $H=[h_1,...,h_\mu]$.
The recovery $\widehat{w}(\omega)=\widehat{w}^H(\omega)$ of $w=Bx$ by the estimate, observation being $\omega$,  is
\beq
\wh{w}(\omega)=B\wh{x}(\omega),\;\;\wh{x}(\omega)\in\Argmin_{x\in X} \|H^T[\omega-Ax]\|_\infty.
\eeq
The error of the estimate $\wh{w}$ allows for the simple upper bound.
\begin{theorem}\label{the1} {\rm \cite[Proposition 5.1]{PUP}} Given $\delta>0$ and a $\delta$-admissible contrast matrix $H\in \bR^{m\times\mu}$, let
\[
\mP[H]=\max\limits_{x} \left\{\|Bx\|_\theta:\, x\in X_\s,\;\|H^TAx\|_\infty\leq1\right\}.
\]
Then for every $x\in X$ the $P_x$-probability of the event
\[
\|\widehat{w}(Ax+\xi_x)-Bx\|_\theta> 2\mP[H]\]
does not exceed $\mu\delta$. Equivalently, for $\epsilon=\mu\delta$ the $\epsilon$-risk
of the polyhedral estimate $\widehat{w}^H(\cdot)$ does not exceed $2\mP[H]$.
\end{theorem}
\section{Design of presumably good polyhedral estimate}\label{secdesign}
\subsection{The strategy}\label{sectStrag}
When designing a polyhedral estimate, the goal is to specify, given $X$, $\epsilon\in(0,1)$, and the family of norms $\pi_\delta$, the contrast matrix $H$ resulting in as small as possible upper bound $\mP[H]$  on the $\epsilon$-risk of the associated with $H$ polyhedral estimate. As stated, this goal is usually unachievable because computing $\mP$ for a given $H$ amounts to maximizing a convex function and typically is a computationally intractable task; in addition, $\mR[H]$ not necessarily is convex in $H$. Instead, we develop  ``presumably good'' designs and our strategy is as follows.
\medskip\par$\bullet$ We start with the following simple
\begin{observation}\label{obs} Let $\theta_*={\theta\over 2-\theta}$. Given $\delta>0$,  let $\zeta\in\bR^\nu$ and $U,S\in\bS^n$ be such that
\begin{equation}\label{LMI}\left[\begin{array}{c|c}U+S&B^T\cr\hline B&\Diag\{\zeta\}\cr\end{array}\right]\succeq0,\;\;
\|\zeta\|_{\theta_*}\leq 1
\end{equation}
Let also $H_1\in\bR^{m\times\mu_1}$, $H_2\in\bR^{m\times \mu_2}$ be $\delta$-admissible contrast matrices. Then
\begin{align}
\mP^2[[H_1,H_2]]&=\max_x\left\{\|Bx\|_\theta^2:x\in X_\s,\|[H_1,H_2]^TAx\|_\infty\leq1\right\}\nn
&\leq\mP[U|\wt{X},H_1]+\mP[S|X_\s,H_2]
\label{observe}
\end{align}
where for $V\in\bS^n$, convex compact set $\Xi\subset\bR^n$ and $G\in\bR^{N\times m}$
\[
\mP[V|\Xi,G]:=\max\limits_{x} \left\{x^TVx: x\in \Xi,\; \|G^TAx\|_\infty\leq1\right\}.
\]
\end{observation}
\medskip\par
$\bullet$ By Theorem \ref{the1}, \rf{observe} implies that the squared $[\mu_1+\mu_2]\delta$-risk of the polyhedral estimate $\widehat{w}_{[H_1,H_2]}$ is at most  $\mP[U|\wt{X},H_1]+\mP[S|X_\s,H_2]$.
\medskip
\par
$\bullet$  We derive upper bounds  on $\mP[U|\wt X,H_1]$ and ${\mP}[S|X_\s,H_2]$ allowing for computationally efficient processing of the problem of optimizing the sum of these bounds in $U$, $S$ $\zeta$ satisfying (\ref{LMI}) and in $\delta$-admissible contrast matrices $H_1$ and $H_2$, resulting in ``presumably good'' contrast matrix $H=[H_1,H_2]$ and associated polyhedral estimate.\footnote{One should not be surprised by the fact that  $X_\s$ is replaced with a larger set $\wt{X}$ in the first term of the bound (\ref{observe}). The reason is that in the approach we are about to present, when upper bounding the first term of the expression, we only utilize the fact that $X_\s\subset\wt{X}$.}

\subsection{Bounding $\mP[U|\cX,H]$}
Upper-bounding of $\mP[U|\cX,H]$, where $\cX$ is the ellitope (\ref{cX}) and $H$ is an $m\times\mu$ matrix is based on the following observation (cf. \cite[Section 5.1.5]{PUP}):
\begin{proposition}\label{prop1} Let  $W\in\bS^n$, $\gamma\in\bR^K_+$, $\Theta\in\bS^m_+$, $\delta$-admissible $m\times\mu$ matrix $H$ and vector $\lambda\in\bR^m_+$ satisfy the relations
\begin{equation}\label{rel1}
P^TWP\preceq P^TA^T\Theta AP+\sum_k\gamma_kT_k,\;\;
H\Diag\{\lambda\}H^T\succeq\Theta
\end{equation}
(here $P$ and $T_k$ come from (\ref{cX})).
Then
\beq\label{res1}
\mP[W|\cX,H]\leq \phi_\cT(\gamma)+\sum_i\lambda_i
\eeq
where from now on
$$
\phi_Z(g)=\max_{z\in Z}g^Tz
$$
is the support function of a set $Z\subset\bR^N$ and $\cT$ is as in (\ref{cX}).
\end{proposition}
To build a mechanism for a convenient for our ultimate purpose upper-bounding $\mP[W|\cX,H]$, we augment Proposition \ref{prop1} by the following observation (which is a far-reaching extension of \cite[Lemma 5.6]{PUP}).
\begin{proposition}\label{prop2} Let $\varrho(\cdot)$ be a norm on $\bR^M$ with the unit ball---a {\em basic ellitope}
$$
B^{\varrho}:=\{g:\,\varrho(g)\leq1\}=\{g:\,\exists s\in\cS: \,g^TS_\ell g\leq s_\ell,\,\ell\leq L\}
$$
where $\cS$ and $S_k$ are as in the definition of an ellitope.
Let us specify the closed convex cone $\bM\subset\bS^{M}_+\times\bR_+$ as
\begin{align}\label{bH}
\bM&=\{(\Xi,\rho)\in\bS^{M}_+\times\bR^+:\exists s\in\cS:
\Tr(\Xi S_\ell)\leq (\rho/\varkappa) s_\ell,\ell\leq L\}\\
\varkappa&=
2\sqrt{2}\ln(4M^2L).\nonumber
\end{align}
Then
\item[(i)] whenever $\Xi=\sum_j\lambda_jg_jg_j^T$ with $\lambda_j\geq0$ and $\varrho(g_j)\leq1$ $\forall i$, we have
$$
\left(\Xi,\,\varkappa{\sum}_j\lambda_j\right)\in \bM,
$$
\item[(ii)] and ``nearly'' vice versa: when $(\Xi,\rho)\in\bM$, there exist (and can be found efficiently by a randomized algorithm)  $\lambda_j\geq0$ and $g_j$, $j\leq M$, such that
$$
\Xi=\sum_j\lambda_jg_jg_j^T\;\;\mathrm{with}\;\;\sum_j\lambda_j\leq\rho\;\;\mathrm{and}\;\;\varrho(g_j)\leq 1,\;j\leq M.
$$
\end{proposition}

When applying Proposition \ref{prop2} with  $\varrho$ having the basic ellitope
\[
\left\{z:\,\exists s\in\cS^\delta:\, z^TS^\delta_\ell z\leq s_\ell,\,\ell\leq L\right\}
\] as the unit ball, we arrive at the following
\begin{corollary}\label{cor01} Given $\delta>0$, consider the set
\[
\bH_\delta=\{(\Theta,\rho):\exists\, \Xi:(\Xi,\rho)\in\bM_\delta,\;\Theta=S_\delta^T \Xi S_\delta\}
\]
where
\[
\bM_\delta=\{(\Xi,\rho)\in\bS^{M}_+\times\bR^+:\exists s\in\cS^\delta:
\Tr(\Xi S^\delta_\ell)\leq (\rho/\varkappa) s_\ell,\ell\leq L\}
\]
(cf. {\rm(\ref{Bpi})}).
This set is a closed convex cone with the following properties:
\item[(i)] when $\Theta=\sum_j\lambda_jh_jh_j^T\in\bS^m_+$ with $\lambda_j\geq0$ and $\pi_\delta(h_j)\leq 1$ for all $j$, one has
$$
(\Theta,\varkappa\sum_j\lambda_j)\in\bH_\delta;
$$
\item[(ii)] whenever $(\Theta,\rho)\in\bH_\delta$, there exist and can be found by an efficient randomized algorithm
$h_j\in\bR^m$ and $\lambda_j\geq0$, $j\leq M$, such that
\beq
\Theta=\sum_j\lambda_jh_jh_j^T\;\;\mathrm{with}\;\; \sum_j\lambda_j\leq \rho\;\;\mathrm{and}\;\;\pi_\delta(h_j)\leq1,\,j\leq M.
\eeq
\end{corollary}

Let now $(\Theta,\rho)\in (\Theta,\rho)\in\bH_\delta $. Corollary \ref{cor01}.{\em ii}
allows us to convert such $\Theta$, in a computationally efficient fashion, into $\delta$-admissible $m\times M$ matrix $H$ and vector $\lambda\in\bR^M_+$ with $\sum_j\lambda_j\leq\rho$ such that $\Theta=H\Diag\{\lambda\}H^T$. Now, applying  Proposition \ref{prop1} we arrive at the following result.
 \begin{theorem}\label{the2}
 Given $\delta>0$, let  $U\in\bS^n$, $\gamma\in\bR^K_+$, $(\Theta,\rho)\in\bH_\delta$ satisfy the relations
\[
P^TUP\preceq P^TA^T\Theta AP+\sum_k\gamma_kT_k.
\]
One can convert, in a computationally efficient way, the above quantities into a $\delta$-admissible $m\times M$ matrix $H$ such that
\[
\mP[U|\cX,H]\leq \phi_\cT(\gamma)+\rho.
\]
 \end{theorem}
\paragraph{Remark} On inspection, the above reasoning establishes the following fact which seems to be important in its own right:
\begin{proposition} \label{prop:randomize}Let
$$
\cH=\{h\in\bR^m:\exists g\in\bR^M,s\in\cS:h=Sg,g^TS_\ell g\leq s_\ell,\ell\leq L\}
$$
be an ellitope. Consider the sets
\begin{align*}
\bH&=\Conv\{hh^T:h\in\cH\},\\
\bH^+&=\{H:\exists G\in\bS^M_+,s\in\cS:H=S^TGS,\Tr(S_\ell G)\leq s_\ell,\ell\leq L\}.
\end{align*}
Then
\beq\label{then1}
\bH\subset\bH^+\subset \varkappa\bH,\,\,\varkappa=2\sqrt{2}\ln(4M^2L).
\eeq
Furthermore, given $H\in\bH^+$, we can find efficiently by a randomized algorithm $h_j\in\cH$ and $\lambda_j\geq0$, $j\leq  M$, such that
\beq\label{moreover}
H=\sum_j\lambda_jh_jh_j^T\;\;\mathrm{and}\;\; \sum_j\lambda_i\leq \varkappa.
\eeq
\end{proposition}
Note that (\ref{then1}), with independent of $M$ constant $3\ln(\sqrt{3}L)$ in the role of $\varkappa$, is an immediate consequence of the result of \cite[Proposition 4.6]{PUP} which characterizes the quality of semidefinite relaxation of the problem of maximizing quadratic form over an ellitope\footnote{In simple cases, e.g. when $\cH$ is the unit $\|\cdot\|_r$-ball with $2\leq r\leq\infty$, $\varkappa$ can be made just an absolute constant \cite{NesSDP,YuLp}}. The novelty, as compared to that result, is the technique for efficient recovery of representation (\ref{moreover}) described in the proof of Proposition \ref{prop2}.

\subsection{Bounding $\mP[S|X_\s,H]$}\label{sectoX}
Suppose that we are in the situation of Section \ref{sitgoal} and are given $\delta>0$ and  matrix $S\in\bS^n$.
\par
Consider the construction as follows. Let us set
\begin{equation}\label{V}
V=[v_1,...,v_J]\in\bR^{p\times J}
\end{equation}
and
\begin{equation}\label{oY}
\begin{array}{ll}
Y&=\left\{y\in\Conv\{\pm v_1,...,\pm v_J\}:\exists w:y=Qw,Rw\in X_\s\right\}\\
&=\left\{y:\exists \lambda,w:y=V\lambda,\|\lambda\|_1\leq 1,y=Qw,Rw\in X_\s\right\}\\
&\subset\overline{Y}:=\left\{y:\exists \lambda,w,z,t\in\cT:
\begin{array}{l}y=V\lambda,\|\lambda\|_1\leq 1,
y=Qw\\
Rw=Pz,z^TT_kz\leq t_k,k\leq K\\
\end{array}\right\}.\\
\end{array}
\end{equation}
Denoting by $Q^+$ the pseudoinverse of $Q$ and taking into account that $\Ker\,Q=\{0\}$, the relation $y=Qw$ implies that $w=Q^+ y$. Thus, if $y\in Y$, so that $y=Qw$ with $Rw\in X_\s$, then $x:=RQ^+ y$. Vice versa, if $x\in X_\s$, then $x\in\overline{X}$, which, by (\ref{ovX}), implies that $x=Rw$ for $w$ such that $y:=Qw\in\Conv\{\pm v_1,...,\pm v_J\}$, so that $y\in Y$. We conclude that
\begin{equation}\label{conclude}
X_\s=\{x:\exists y\in Y: x=RQ^+ y\}.
\end{equation}
Thus, given $x\in X_\s$, we can find $y\in Y$ such that $x=RQ^+ y$; since $y\in Y$, there exists $\lambda$ with $\|\lambda\|_1\leq1$ such that $y=V\lambda$. Therefore
$$
\begin{array}{l}
x^TSx=y^T\overbrace{[Q^+]^TR^TSRQ^+}^{\overline{S}[S]}y=[\sum_j\lambda_jv_j]^T\overline{S}[S]y.
\leq \max_j|v_j^T\overline{S}[S]y|
\end{array}
$$
Thus,
\begin{equation}\label{thus}
\begin{array}{c}
\forall x\in X_\s\exists y\in Y\subset\overline{Y}: x=RQ^+ y\ \&\ x^TSx\leq \max\limits_{j\leq J}|v_J^T\overline{S}[S]y|,\\
\overline{S}[S]=[Q^+]^TR^TSRQ^+.\\
\end{array}
\end{equation}
Now imagine that for every $j\leq J$ we have at our disposal  $\delta$-admissible matrix $H_j$. Setting $H=[H_1,...,H_J]$, (\ref{thus}) says that
\begin{equation}\label{(*)}
\mP[S|X_\s,H]\leq \max\limits_{j\leq J}\left[ \mq_j(H_j):=\max_y\left[|v_J^T\overline{S}[S]y|: y\in Y,\|H_j^TARQ^+ y\|_\infty\leq 1\right]\right].
\end{equation}
Our approach to optimizing (an upper bound on) $\mP[S|X_\s,H]$ over $\delta$-admissible contrast matrices $H$ is based on the specific policy for optimizing quantities $\mq_j(H_j)$ over $\delta$-admissible matrices $H_j$, summarized in \cite[Observation 5.1]{PUP} and setting $H=[H_1,H_2,...,H_J]$ (at optimum, matrices $H_j$ turn out to be column vectors). In the rest of this section we describe and justify this policy.
\par
Observe that $Y$ is a nonempty convex compact set symmetric w.r.t. the origin. Consider optimization problems
$$
\begin{array}{rcl}
\varsigma_j&=&\min_{g}\left\{\phi_j(S,{g}):=\max_{y\in Y}[v_j^T\overline{S}[S]-{g}^TARQ^+]y+\pi_\delta({g})\right\}\\
&=&\min_{g}\left\{\phi_j(S,{g}):=\max_{y\in Y}|[v_j^T\overline{S}[S]-{g}^TARQ^+]y|+\pi_\delta({g})\right\},\\
\end{array}  \eqno{(P_j)}
$$
where the second equality is due to the symmetry of $Y$ w.r.t. the origin.
Problems ($P_j$) clearly are solvable; let ${g}_j$ be their optimal solutions, and ${h}_j$, $\pi_\delta({h}_j)=1$, be such that $\anc{h}{g}_j=\pi_\delta({g}_j){h}_j$.
 Note that $\phi_j(S,{g})$ is a convex function of $(S,{g})$,
whence $\varsigma_j$ is convex function of $S$. Let
$$
x\in X_\s\ \&\ |{h}_j^TAx|\leq1.
$$
Then for every $y\in Y$ such that $x=RQ^+ y$, setting $S'=\overline{S}[S]$, we have
\begin{align*}
|v_j^TS'y|&\leq |v_j^TS'y-g_j^TARQ^+ y|+|g_j^TARQ^+ y|\\
&\leq \varsigma_j-\pi_\delta(g_j)+|g_j^TA\underbrace{RQ^+ y}_{=x}|=\varsigma_j-\pi_\delta(g_j)+\pi_\delta(g_j)|{h}_j^TAx|\leq \varsigma_j.
\end{align*}
Thus, if $x\in X_\s$ and
$|{h}_j^TAx|\leq1,1\leq j\leq J$, then for every $y\in Y$ such that $x=RQ^+y$ one has
$$
 \left\|[v_1^TS'y;...;v_J^TS'y]\right\|_\infty\leq\varsigma:=\max_j\varsigma_j.
$$
This combines with (\ref{thus}) to imply that whenever $x\in X_\s$ is such that $|{h}_j^TAx|\leq1$, $j\leq J$, one has $x^TSx\leq \varsigma,$
or, which is the same,
\[
\mP[S|X_\s,[{h}_1;...;{h}_J]]\leq \varsigma.
\]

\subsection{Putting things together}\label{sec:together}
We have described computationally efficient techniques for upper bounding $\mP(\cdot|\cX,\cdot)$ and $\mP(\cdot|X_\s,\cdot)$ and thus are now ready to implement the strategy for building presumably good polyhedral estimate outlined in Section \ref{sectStrag}.

Given $\delta>0$, consider the convex optimization problem

\begin{equation}\label{fprob}
\begin{array}{l}
\min\limits_{\Theta,\zeta,\rho,\gamma,\atop
\varsigma,U,S,\{{g}_j,j\leq J\}}
\bigg\{\phi_\cT(\gamma)+\rho+\varsigma:\\
\qquad\left.\begin{array}{l}
\left[\begin{array}{c|c}U+S&B^T\cr\hline B&\Diag\{\zeta\}\cr\end{array}\right]\succeq0,\|\zeta\|_{\theta_*}\leq1,\\
(\Theta,\rho)\in\bH_\delta,\gamma\geq0,P^TUP\preceq P^TA^T\Theta AP+\sum_k\gamma_kT_k,\\
\max_{y\in Y}\left[v_j^T[Q^+]^TR^TSRQ^+-{g}_j^TARQ^+\right]y+\pi_\delta({g}_j)\leq\varsigma,\,j\leq J,\\
\end{array}\right\}
\end{array}
\end{equation}
and let $\Theta,...,\{{g}_j,\,j\leq J\}$ be its feasible solution. By Theorem \ref{the2}, we can convert, into a computationally efficient fashion, component $\Theta$ of this feasible solution into
$\delta$-admissible $m\times  M$ matrix $H_1$ such that
\[
\mP[U|\cX,H_1]\leq \phi_\cT(\gamma)+\rho.
\]
Next, as we have just seen in Section \ref{sectoX}, when converting ${g}_j$, $j\leq J$, into a $\delta$-admissible $m\times J$ matrix $H_2$ we get
\[
\mP[S|X_\s,H_2]\leq \sqrt{\varsigma}.
\]
When applying Observation \ref{obs} we conclude that setting $H=[H_1,H_2]$, we obtain a $\delta$-admissible $m\times[M+J]$-matrix $H$ such that
\beq\label{synthesis}
\mP^2[H]\leq \phi_\cT(\gamma)+\rho+\varsigma,
\eeq
implying by Theorem \ref{the1} that for $\epsilon=\delta[M+J]$
\beq\label{final}
\Risk_{\epsilon}[\widehat{w}^H|X]\leq 2\sqrt{\phi_\cT(\gamma)+\rho+\varsigma},
\eeq
Our last step is to get rid of $\max_{y\in Y}$ in (\ref{fprob}). To this end we first replace it with $\max_{y\in\overline{Y}}$, thus reducing the feasible set of (\ref{fprob}) (recall that $Y\subset\overline{Y}$) and therefore preserving the validity of (\ref{synthesis}) and (\ref{final}).
Setting
$
\bT=\cl\{[z;s]:s>0,z/s\in\cT\},
$
the dual to $\bT$ cone is
$
\bT_*=\{[\rho;\sigma]:\sigma\geq\phi_\cT(-\rho)\}.
$ Let us put
\begin{equation}\label{Ek}
E_k=2T_k^{1/2},%
\end{equation}
so that the inequality $z^TT_kz\leq t_k$ is equivalent to conic inequality $[E_kz;t_k-1;t_k+1]\in\bL^{N+2}$, where $\bL^r$ is the $r$-dimensional Lorentz cone.
Thus, for $d\in \bR^p$,
\[
\max\limits_{y\in\overline{Y}}d^Ty=\max\limits_{y,\lambda,w,z,t}\left\{d^Ty:
\begin{array}{l}
y-V\lambda=0, \,y-Qw=0,Rw-Pz=0,\,
\|\lambda\|_1\leq 1,\\~
[E_kz;t_k-1;t_k+1]\in\bL^{N+2},k\leq K,[t;1]\in\bT,
\end{array}\right\}
\]
and by the Conic Duality Theorem
\begin{align*}
\max\limits_{y\in\overline{Y}}d^Ty
&=\min\limits_{{\{\alpha,\beta\}\subset\bR^p,\eta\in\bR^n,
\delta\in\bR^J,\atop[\epsilon_k;\phi_k;\psi_k]\in\bL^{N+2},k\leq K,[\rho;\sigma]\in\bT_*}}
\Big\{\|\delta\|_\infty +\sum_k[\psi_k-\phi_k]+\sigma:\\&\left.\qquad\qquad\qquad\qquad\qquad\qquad
\begin{array}{l}\alpha+\beta={d},\,\delta-V^T\alpha=0,\,R^T\eta-Q^T\beta=0,\\
\sum_kE_k^T\epsilon_k+P^T\eta=0,\;
\phi+\psi=-\rho
\end{array}\right\}\\
&=\min\limits_{{\beta\subset\bR^p,\eta\in\bR^n,
\atop\{\epsilon_k\in\bR^N,k\leq K\},\{\phi,\psi,\rho\}\subset\bR^K}}\Big\{\|V^T[{d}-\beta]\|_\infty +\sum_{k=1}^K[\psi_k-\phi_k]+
\phi_\cT(\phi+\psi):\\
&\left.\qquad\qquad\qquad\qquad\qquad\qquad\begin{array}{l}R^T\eta-Q^T\beta=0,\;\sum_kE_k^T\epsilon_k+P^T\eta=0\\
\psi_k\geq \|[\epsilon_k;\phi_k]\|_2,\,k\leq K.
\end{array}\right\}
\end{align*}
We arrive at the following summary of the situation.
\par
 {\em In the situation of Section \ref{sitgoal} and given $\delta>0$, consider convex optimization problem
\color{black}
\begin{equation}\label{fprobe}
  \begin{array}{l}
\Opt=\min\limits_{{\Theta,\zeta,\zeta,\rho,\gamma,
\varsigma,U,S\{{g}_j,j\leq J\},\atop
\{\beta^j,\eta^j,\epsilon^j_k,k\leq K,\phi^j,\psi^j\}_{J\leq J}}}\Big\{\phi_\cT(\gamma)+\rho+\varsigma:
\\
\quad
\left.\begin{array}{l}
\left[\begin{array}{c|c}U+S&B^T\cr\hline B&\Diag\{\zeta\}\cr\end{array}\right]\succeq0,\|\zeta\|_{\theta_*}\leq1,\\
(\Theta,\rho)\in\bH_\delta,\gamma\geq0,\\P^TUP\preceq P^TA^T\Theta AP+\sum_k\gamma_kT_k,\\
\|V^T\left[[Q^+]^TR^TSRQ^+v_j-[Q^+]^TR^TA^T{g}_j-\beta^j\right]\|_\infty\\
\qquad+\sum_k[\psi^j_k-\phi^j_k]+\phi_\cT(\phi^j+\psi^j)+\pi_\delta({g}_j)
\leq\varsigma\\
R^T\eta^j-Q^T\beta^j=0,\sum_kE_k^T\epsilon^j_k+P^T\eta^j=0,j\leq J\\
\psi^j_k\geq \|[\epsilon^j_k;\phi^j_k]\|_2,\,k\leq K,\,j\leq J.
\end{array}\right\}
\end{array}
\end{equation}
(for notation, see Section \ref{sitgoal}, Corollary \ref{cor01},  (\ref{V}), (\ref{oY}), and (\ref{Ek})).

A feasible solution $\Theta,... ,$ to this problem can be efficiently converted into $m\times[J+M]$ $\delta$-admissible contrast matrix $H$ such that with $\epsilon=[M+J]\delta$, the $\epsilon$-risk of the associated with $H$ polyhedral estimate $\widehat{w}^H$ satisfies
$$
\Risk_{\epsilon}[\widehat{w}^H|X]\leq \sqrt{\phi_\cT(\gamma)+\rho+\varsigma}.
$$
} 
\paragraph{Remark} We have assumed that the recovery error is measured in $\|\cdot\|_\theta$, $\theta\in[1,2]$. This assumption may be relaxed---$\|\cdot\|_\theta$ can be replaced with any norm $\|\cdot\|$ admitting the following description:
$$
\|u\|^2=\inf_{\Delta\in\cD,\,\Delta\succ0}\langle u,\Delta^{-1}u\rangle
$$
where $\cD\subset\bS^{\nu}_+$ is a convex compact set containing a positive definite matrix.
The above summary remains valid when variable  $\zeta\in \bR^\nu$ in (\ref{fprobe}) is replaced with $\Delta\in\bS^\nu$ and constraint $\left[\begin{array}{c|c}U+S&B^T\cr\hline B&\Diag\{\zeta\}\cr\end{array}\right]\succeq0,\|\zeta\|_{\theta_*}\leq1$,---with \[\left[\begin{array}{c|c}U+S&B^T\cr\hline B&\Delta\cr\end{array}\right]\succeq0,\;\Delta\in\cD.
\]
In the situation we have considered so far, $u\in\bR^\nu$, $\langle\cdot,\cdot\rangle$ is the standard inner product on $\bR^\nu$, and $\|\cdot\|$ is $\|\cdot\|_\theta$ with $\theta\in[1,2]$, that is, $\cD=\{\Diag\{\delta\}:\delta\geq0,\|\delta\|_{\theta_*}\leq1\}$, $\theta_*={\theta\over2-\theta}$. In somewhat ``more exotic'' example one may consider,  linear mapping $x\mapsto Bx$ takes values in $\bR^{\nu\times \mu}$, $\langle\cdot,\cdot\rangle$ is the Frobenius inner product on $\bR^{\nu\times \mu}$, and $\|\cdot\|$ is the Schatten $|\cdot|_\theta$-norm on $\bR^{\nu\times\mu}$ with $\theta\in[1,2]$ (i.e.,  $|z|_\theta=\|\sigma(z)\|_\theta$ where $\sigma(z)$ is the singular spectrum of  $z\in \bR^{\nu\times\mu}$), corresponding to $\cD=\{\Delta\in\bS^\nu_+:\;|\Delta|_{\theta_*}\leq1\}$.
\section{Numerical illustration}\label{sec:num}
In this section, to illustrate the numerical performance of approach to polyhedral estimate design described in Section \ref{secdesign}, we present results of two ``proof of concept'' simulation experiments.
 \par
 The setting of our {\em first experiment} is as follows. We consider the observation model
 \[
 \omega=Ax+\xi_x\in\bR^m,
 \]
(cf. \rf{eq:mod1}) in which
\begin{itemize}
 \item unknown signal $x$ is known to belong to
 \[X=X_\s=\{x\in\bR^n:\,\|x\|_1\leq1,\|x\|_2\leq \rho_2,\|x\|_\infty\leq \rho_\infty\},
  \]an intersection of the ellitope
 \begin{align*}
 \cX&=\{x:\|x\|_2\leq\rho_2,\|x\|_\infty\leq\rho_\infty\}\\
 &=\{x\in\bR^n:\rho_\infty^{-2}x_k^2\leq1,k\leq n,\rho_2^{-2}x^Tx\leq
 1\}\\&\qquad [K=n+1,\;\cT=\{t:0\leq t_k\leq 1,k\leq K\}]
 \end{align*}
 and the $\ell_1$-ball
 \[
 \oX=\{x:\|x\|_1\leq\rho_1\}\ [R=Q=I_n,J=n,[v_1,...,v_J]=\rho_\infty I_n]\\
\]
 \item $m=n$, $\xi_x\sim\cN(0,\sigma^2I_n)$, whence
 $$
 \Prob\{|h^T\xi_x|>1\}\leq\delta\in(0,1) \;\Leftrightarrow \;\sigma\chi_\delta\|h\|_2\leq 1
 $$
 where $\chi_\delta$ is the $(1-\delta/2)$-quantile of the standard normal distribution. We put
 $$
\pi_\delta(h)=\sigma\chi_\delta\|h\|_2.
 $$
 \item Recovery error is measured in $\|\cdot\|_2$, that is, $\theta=2$.
 \end{itemize}
 In the present setting (\ref{fprobe}) becomes the convex optimization problem
\begin{align}\label{fprobeG}
\Opt&=\min\limits_{{\Theta,\gamma,
\varsigma,\atop
U,S\{\alpha_j,\beta_j,{g}_j,j\leq n\}}}\left\{\sum_{k=1}^{n+1}\gamma_k+\sigma^2\chi_\delta^2\Tr(\Theta)+\varsigma:\right.\\
&\left.\qquad\begin{array}{l}
U+S\succeq B^TB,\;
\Theta\succeq0,\gamma\geq0,\\U\preceq A^T\Theta A+\rho_2^{-2}\gamma_1I_n+\rho_\infty^{-2}\Diag\{\gamma_2,\gamma_2,...,\gamma_{n+1}\}\\
\rho_1\|\rho_1\Col_j[S]-A^T{g}_j-\alpha_j-\beta_j\|_\infty+\rho_\infty\|\alpha_j\|_1\\\qquad\qquad+\rho_2\|\beta_j\|_2+\sigma\chi_\delta\|{g}_j\|_2\leq
\varsigma,\;j\leq n.
\end{array}\right\}\nonumber
\end{align}
In this case, to  convert $\Theta$ into $\delta$-admissible contrast matrix $H$ it suffices to compute the eigenvalue decomposition of $\Theta$ and then set the columns of $H$ to scaled to $\pi_\delta$-norm  equal to 1 eigenvectors of $\Theta$.

\begin{figure}[!htb]
\begin{center}
\includegraphics[width=0.6\textwidth]{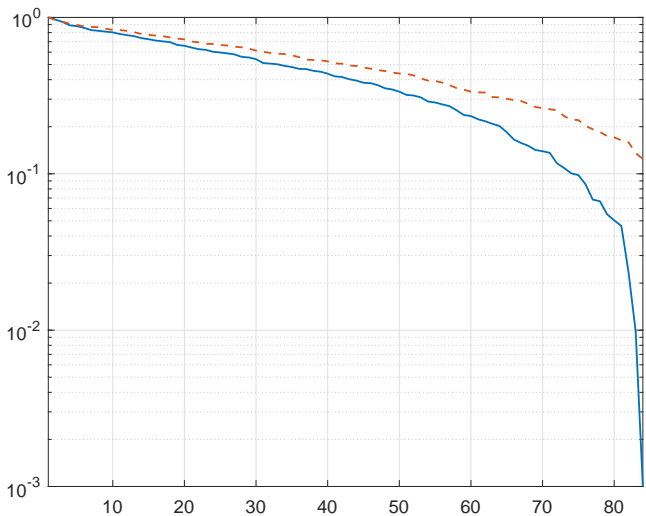}
\end{center}
\caption{\label{fig:0} Singular spectrum  of $A$ (solid line) and $B$ (dashed line).}
\end{figure}

\begin{figure}[!htb]
\begin{tabular}{cc}
\includegraphics[width=0.45\textwidth]{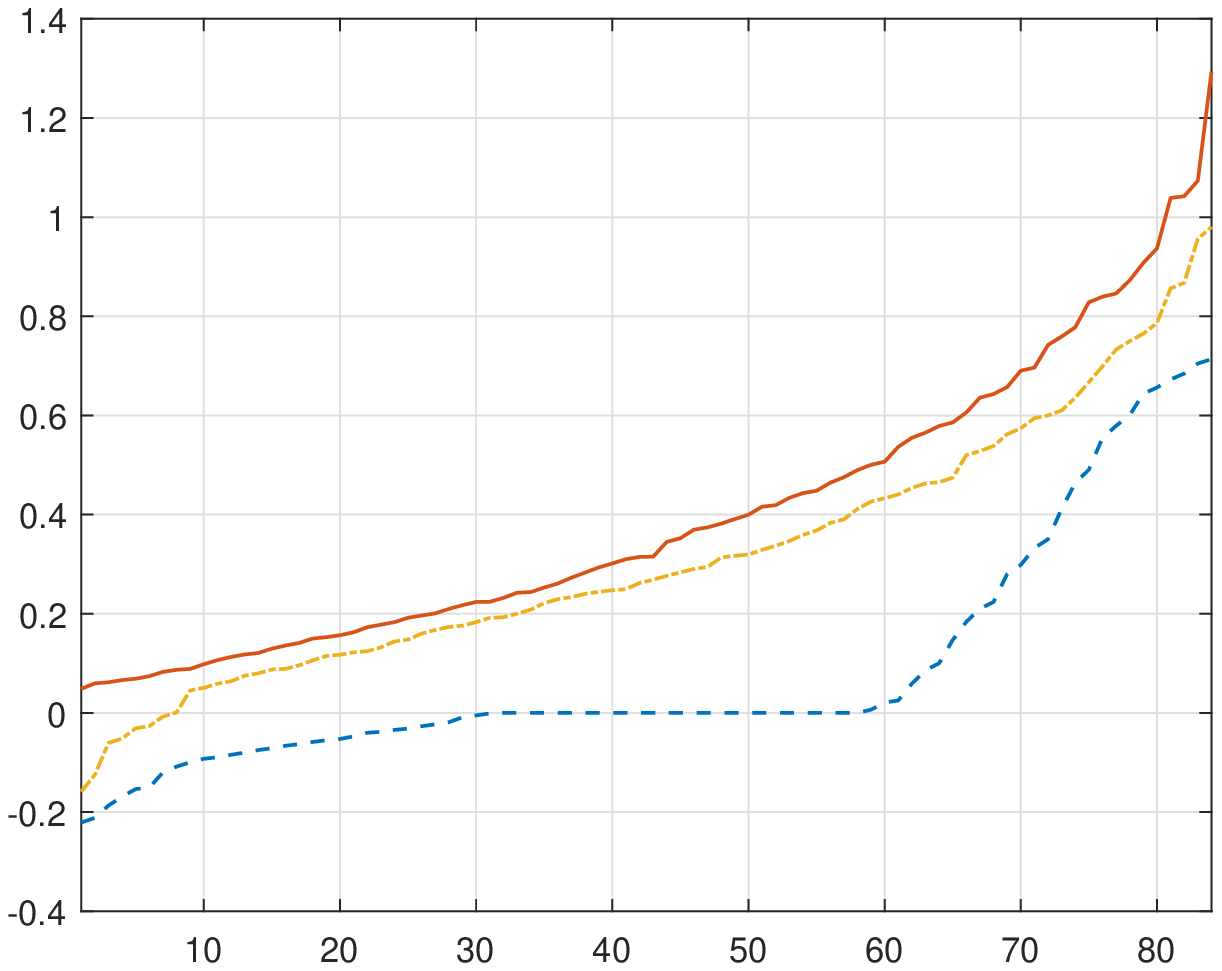}&
\includegraphics[width=0.45\textwidth]{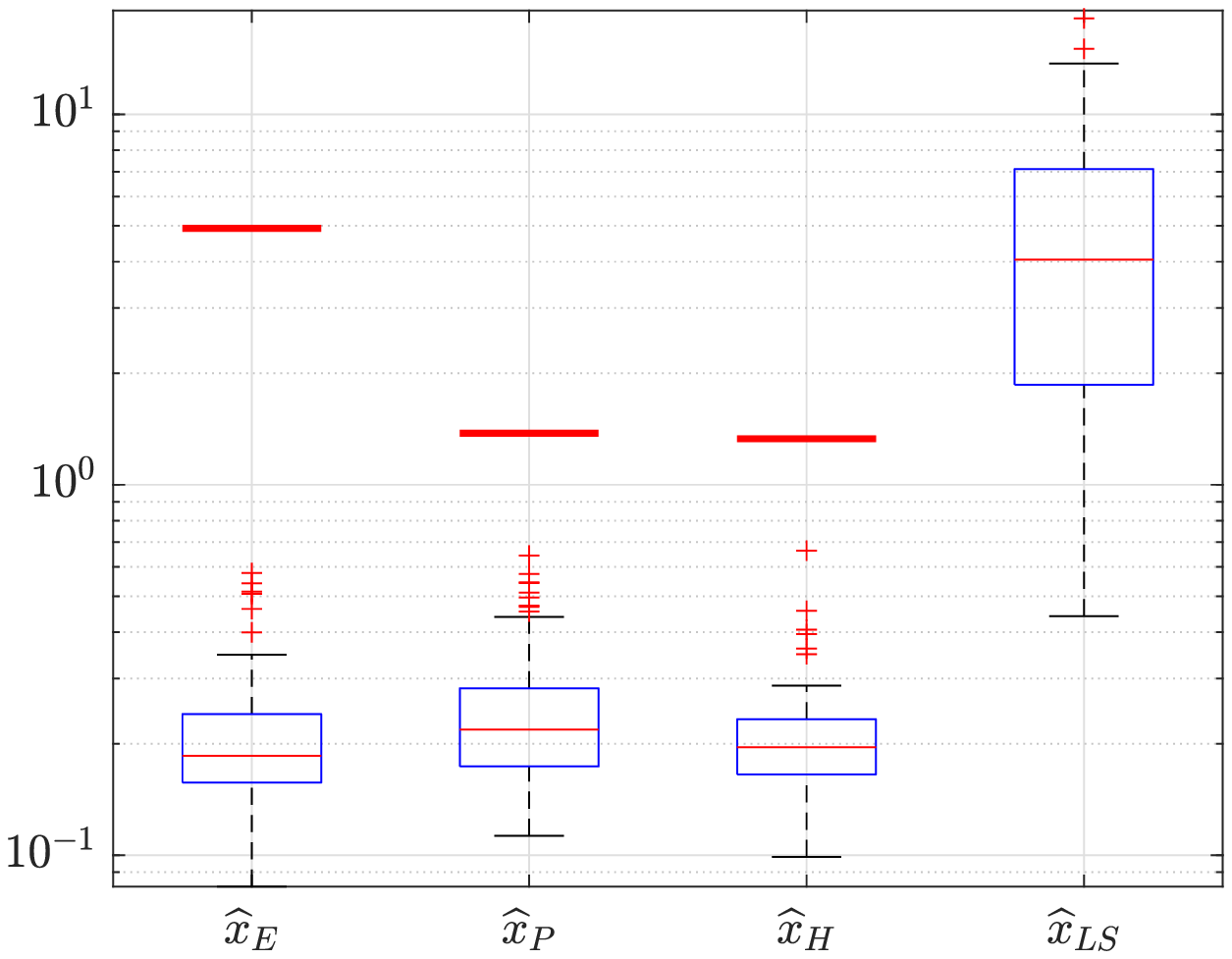}
\end{tabular}
\caption{\label{fig:1} Left plot: eigenvalues of $U$ (dashed line), $S$ (dash-dot line) and $U+S$ (solid line); right plot: distributions of $\|\cdot\|_2$-recovery errors
and theoretical upper bounds on $\Risk_{0.01}$ (red horizontal bars).}
\end{figure}

In the experiment we report on below, $n=m=64$, $\nu=126$, $\rho_1=10$, $\rho_2=8.5$, $\rho_\infty=7$, $\sigma=0.01$, sensing matrix $A$ (condition number 1e3) and matrix $B$ (condition number 8e0) are generated randomly; singular spectra of these matrices are presented in Figure \ref{fig:0}.
We process (\ref{fprobeG}) using {\tt Mosek} commercial solver \cite{mosek} via {\tt CVX} \cite{cvx2014}. Figure \ref{fig:1} illustrates the results of the computation: in the left plot we present ``decomposition'' of $B^TB\preceq I_n$ into matrices $U$ and $S\succeq B^TB-U$. Quite surprisingly, matrices $S$ and $U$ are {\em not} positive semidefinite matrices. In each experiment,  we compute $N=100$ recoveries
$\wh x_H$ of randomly selected signals $x\in X$ with matrix $H$ derived from the optimal solution to \rf{fprobeG} along with two estimates
which only utilize partial information about $x$: when computing estimate $\wh x_E$ we only use $x\in \cX$ (we set $S=0$ in this case), while recovery $\wh x_P$ only utilizes $x\in \overline X$ (we set $U=0$).
The results are presented in the right plot of Figure \ref{fig:1}; for comparison, we provide the error boxplot of the Least Squares estimate $\wh x_{LS}=A^{-1}\omega $.
\par
\begin{figure}[!htb]
\begin{center}
\includegraphics[width=0.6\textwidth]{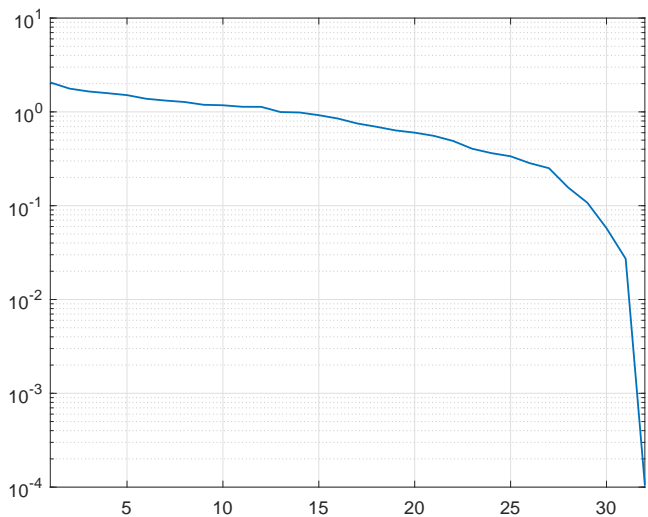}
\end{center}
\caption{\label{fig:2} Singular spectrum  of $A$.}
\end{figure}
\begin{figure}[!htb]
\begin{tabular}{cc}
\includegraphics[width=0.45\textwidth]{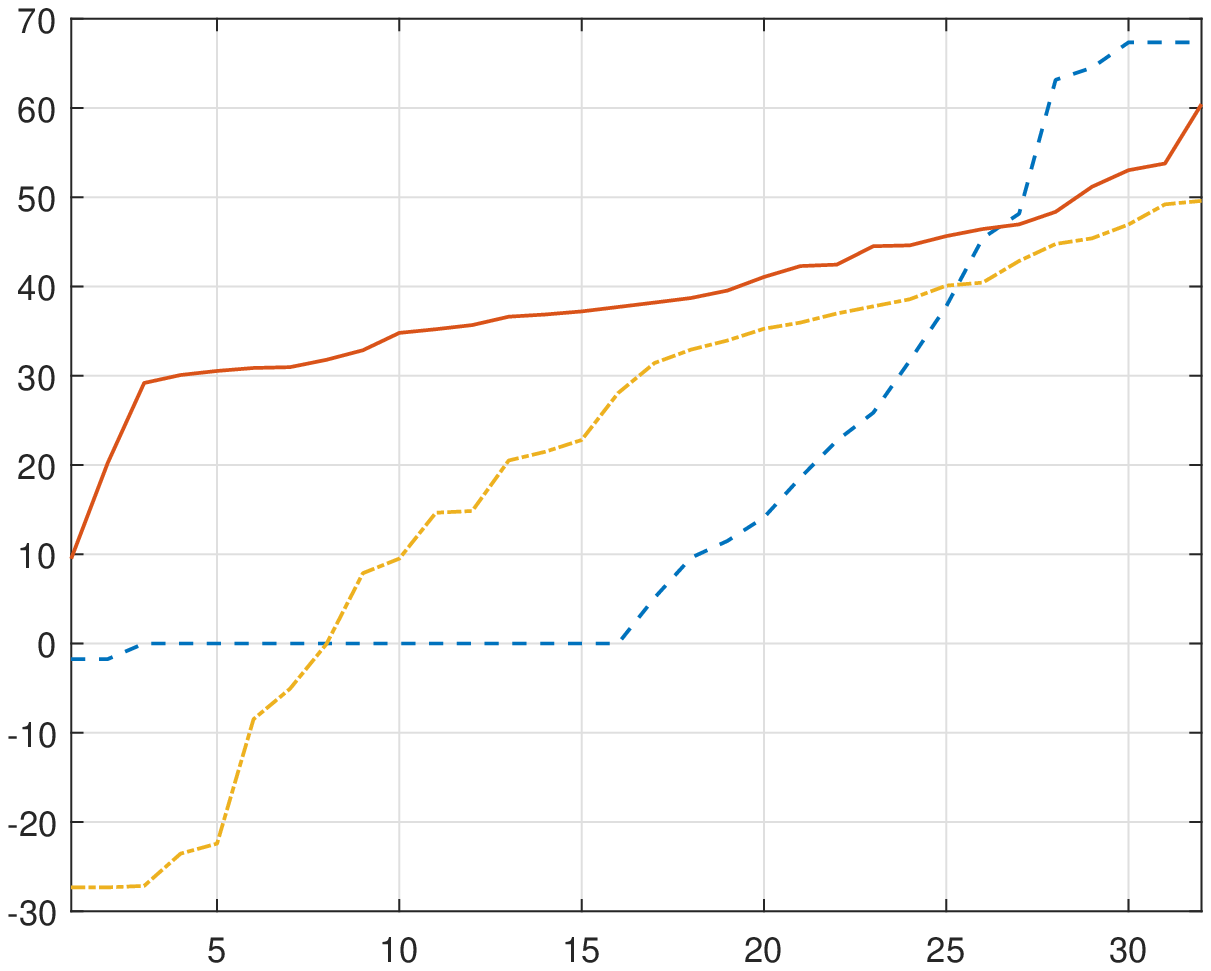}&
\includegraphics[width=0.45\textwidth]{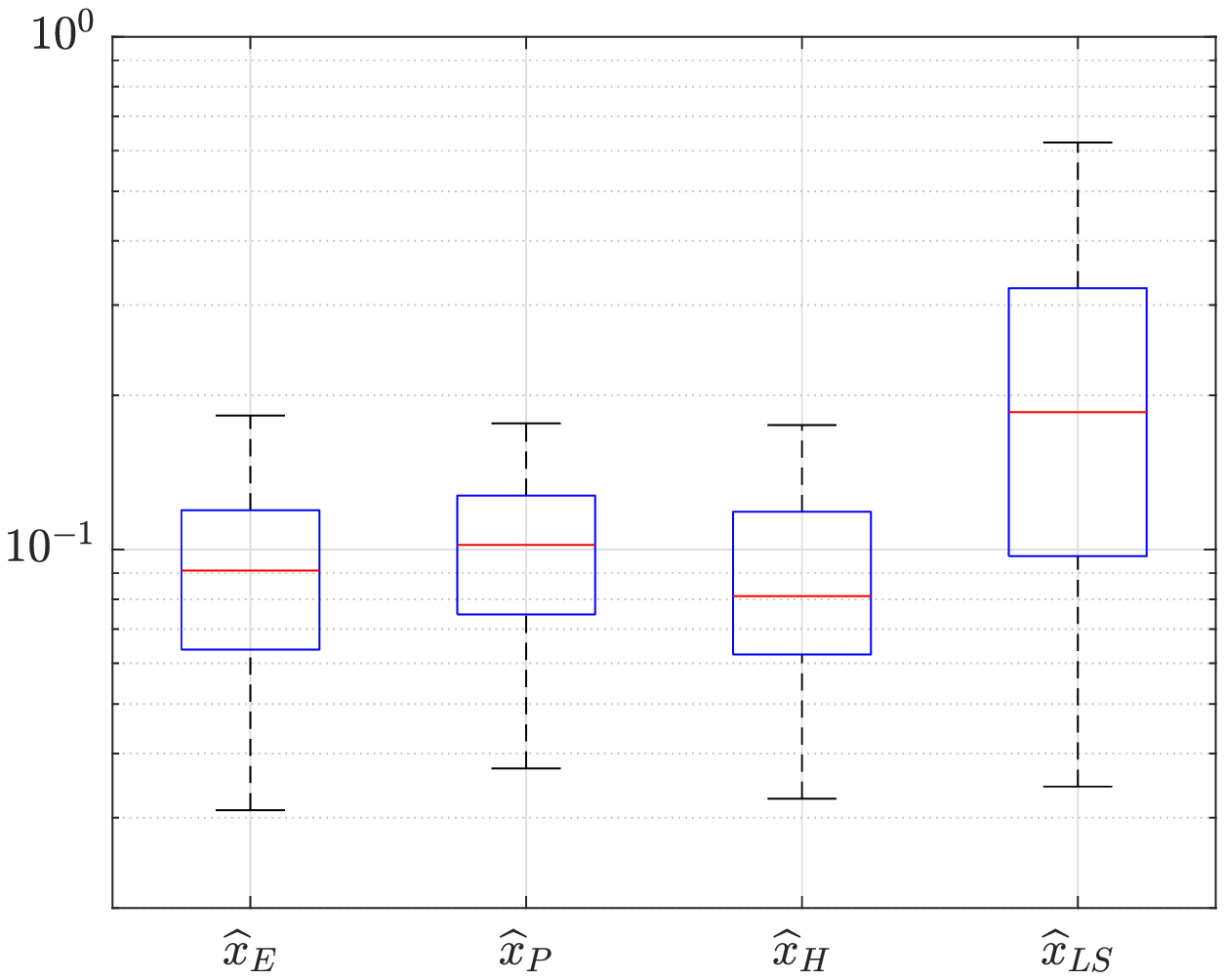}
\end{tabular}
\caption{\label{fig:3} Left plot: eigenvalues of $U$ (dashed line), $S$ (dash-dot line) and $U+S$ (solid line); right plot: distributions of $\|\cdot\|_1$-recovery errors for randomly sampled $x\in X$. Corresponding theoretical upper bounds on $\Risk_{0.01,1}$ are 5.42 for estimate $\wh x_E$, 1.72 for estimate $\wh x_P$, and $1.61$ for estimate $\wh x_H$.}
\end{figure}

{\em Our second illustration} deals with the mixture-sub-Gaussian regression model described in Section \ref{sec:numeric}. In our experiment, $N=1$e$4$, $n=d=32$, $\theta=1$, $B=I_n$, recovery error is measured in $\ell_1$-norm ($\theta=1$), and the ellitope
\[
\cX=\{x\in\bR^n:\,\|x\|_2\leq1,\|x\|_\infty\leq 1/2\}.
\]
Distribution of $\omega_t\in \bR^n$ is a mixture of normal distributions with parameters $(a_i, \Theta_i)$, $i=1,..., n$, selected at random with  $\|a_i\|_2=1$  and positive semidefinite $\Theta_i$ of unit spectral norm. Figure \ref{fig:2} shows the typical singular spectrum of the ``regression matrix'' $A\in \bR^{n\times n}$.
Same as above, in this experiment we compute three polyhedral estimates corresponding to optimal solutions to \rf{fprobe}: ``full'' solution with nonvanishing matrices $U$ and $S\succeq I_n-U$ (estimate $\wh x_H$), and solutions with zeroed matrices $S$ (estimate $\wh x_E$) an $U$ (estimate $\wh x_P$).
Note that, in the present setting,  converting $\Theta$ into a $\delta$-admissible contrast matrix $H$ requires implementing ``full'' randomized procedure as described in the proof of Proposition \ref{prop2} (cf. Proposition \ref{prop:randomize}).
Typical results are presented in Figure \ref{fig:3}.

\appendix
\section{Proofs}
\subsection{Proof of Observation \ref{obs}}
Indeed, under the premise of the observation, 
let $x$ satisfy $\|[H_1,H_2]^TAx\|_\infty\leq1$ and $x\in X_\s$. As is immediately seen, we have \[\|a\|_\theta^2=\inf_{\gamma}\left\{\sum_i{a_i^2\over \gamma}:\,\gamma>0,\;\|\gamma\|_{\theta_*}\leq1\right\},
\] whence
\begin{align*}
\|x\|_\theta^2&\leq \lim_{\epsilon\to+0}\sum_i[Bx]_i^2/[\zeta_i+\epsilon]=\lim_{\epsilon\to+0}x^TB^T\Diag\{\zeta+\epsilon[1;...;1]\}^{-1}Bx\\
&\leq x^T[U+S]x\qquad\hbox{[by (\ref{LMI}) and the Schur Complement Lemma]}\\
&=x^TUx +x^TSx\leq \mP[U|\wt{X},H_1]+\mP[S|X_\s,H_2]
\end{align*}
due to $x\in X_\s\subset\wt{X}$ and definition of $\mP[V|\Xi,G]$.\qed

\subsection{Proof of Proposition \ref{prop1}.}
Let $x\in\cX$ be such that $\|H^TAx\|_\infty\leq1$ and $h_1,...,h_\mu$ be the columns of $H$. Since $x\in\cX$, we have $x=Pz$ for some $z\in\bR^N$ such that $z^TT_kz\leq t_k$, $k\leq K$, for some $t\in\cT$. By the first relation in (\ref{rel1}) we have the first inequality on the following chain:
\begin{align*}
x^TWx&=z^TP^TWPz\leq z^T\left(P^TA^T\Theta A P+\sum_k\gamma_kT_k\right)z\\&=x^TA^T\Theta Ax+\sum_k\gamma_kz^TT_kz\\
&\leq \sum_i\lambda_ix^TA^Th_ih_i^TAx+\sum_k\gamma_kt_k \hbox{\ [by the second relation in (\ref{rel1})]} \\
&\leq \sum_i\lambda_i+\phi_\cT(\gamma)\hbox{\ [due to $\|H^TAx\|_\infty\leq1$ and $t\in\cT$]}.
\end{align*}
Taking in the resulting inequality the supremum over $x\in\cX$ satisfying $\|H^TAx\|_\infty\leq1$. we get (\ref{res1}). \qed

\subsection{Proof of Proposition \ref{prop2}.}
{\em (i)}: Let $\lambda_j\geq0$, $\varrho(g_j)\leq 1$, $i\leq N$, and $\Xi=\sum_j\lambda_jg_jg_j^T$. Then for every $j$ there exists $s^j\in\cS$ such that $g_j^TS_\ell g_j\leq s^j_\ell$, $\ell\leq L$, so that when $\lambda\neq0$, setting $s=[\sum_j\lambda_j]^{-1}\sum_j\lambda_js_i\in\cS$, we have
$$
\Tr(\Xi S_\ell)\leq [\sum_i\lambda_i]s_\ell,\ell\leq L,
$$
implying that $(\Xi,\varkappa\sum_i\lambda_i)\in \bM$. The latter inclusion is true as well when $\lambda=0$.
\\
{\em (ii)}: Let $(\Xi,\rho)\in\bM$, and let us prove that $\Xi=\sum_{i=1}^{M} \lambda_ig_ig_i^T$ with $\varrho(g_i)\leq 1$, $\lambda_i\geq0$, and $\sum_i\lambda_i\leq\rho$. There is nothing to prove when $\rho=0$, since in this case $\Xi=0$ due to $(\Xi,0)\in\bM$ combined with $\sum_\ell S_\ell\succ0$. Now let $\rho>0$,  $Z=\Xi^{1/2}$, and let $O$ be the orthonormal $M\times M$ matrix of $M$-point Discrete Cosine Transform, so that all entries in $O$ are in magnitude $\leq\kappa:=\sqrt{2}$. Observe that for deterministic $a\in \bR^M$ and Rademacher random vector $\epsilon\in\bR^M$ (i.e., random vector with independent entries taking values $\pm1$ with probability $1/2$) one has
\beq\label{Tolya}
t\geq0\Rightarrow \Prob\left\{\left({\sum}_ia_i\epsilon_i\right)^2>t{\sum}_ia_i^2\right\}\leq 2\exp\{-t/2\}.
\eeq
Now let $S\in\bS^M_+$, and let $\varsigma=\Tr(\Xi S)$. Setting $Z^\epsilon=Z\Diag\{\epsilon\}O$, we have $Z^\epsilon[Z^\epsilon]^T\equiv \Xi$.
Let now
$$
r^\epsilon_j=\Col_j^T[Z^\epsilon]^TS\Col_j[Z^\epsilon].
$$
For $S=E^TDE$ with orthogonal $E$ and $D=\Diag\{\sigma_i\}$ and $F^\epsilon=EZ^\epsilon=F\Diag\{\epsilon\}O$ with $F=EZ$, we have
\beq\label{lab1}
\varsigma=\Tr(SZZ^T)=\Tr(D[EZ][EZ]^T) = \sum_{i,k}\sigma_iF_{ik}^2.
\eeq
On the other hand,
\begin{align*}
r^\epsilon_j&=\Tr(S\Col_j[Z^\epsilon]\Col_j^T[Z^\epsilon])=\Tr(D[E\Col_j[Z^\epsilon]][E\Col_j[Z^\epsilon]]^T)\\&=
\Tr(D\Col_j[F^\epsilon]\Col_j^T[F^\epsilon])
  =\sum_i\sigma_i[F^\epsilon_{ij}]^2
={\sum}_i\sigma_i\left[{\sum}_kF_{ik}\epsilon_kO_{kj}\right]^2.
\end{align*}
By (\ref{Tolya}), for every $t>0$ probability of the event
$$
[\sum_kF_{ik}\epsilon_kO_{kj}]^2\leq t\sum_kF^2_{ik}O^2_{kj}\,\,\forall i,j
$$
is at least
$$
p=1-2M^2\exp\{-t/2\}.
$$
When this event takes place, we have
$$
r^\epsilon_j\leq t{\sum}_i\sigma_i{\sum}_kF^2_{ik}O^2_{kj}\leq {t\kappa\over M}{\sum}_i{\sum}_k\sigma_iF^2_{ik}\leq {t\kappa\over M} \varsigma
$$
(we have used (\ref{lab1})). It follows that for every $t>0$ probability of the event
$$
\Col_j^T[Z^\epsilon]^TS\Col_j[Z^\epsilon]\leq {t\kappa\over M}\varsigma,\,\,j=1,...,m,
$$
is at least $1-2M^2\exp\{-t/2\}$.
\par
Due to $(\Xi,\rho)\in\bM$, we have  $\Tr(\Xi S_\ell)\leq (\rho/\varkappa)s_\ell$, $\ell\leq L$, for some $s\in\cS$, and by the above the probability of the event
$$
\forall j\leq M: \Col_j^T[Z^\epsilon]^TS_\ell \Col_j[Z^\epsilon]\leq {r\rho\kappa\over\varkappa M}s_\ell,  \ell\leq L,
$$
is at least
$$
1-2M^2L\exp\{-r/2\}
$$
When this event takes place, we have for $\lambda_j={\rho/M}$
$$
\Xi={\sum}_{j=1}^M\lambda_jg_jg_j^T,\,\,g_j=\sqrt{M/\rho}\Col_j[Z^\epsilon],
$$
with $\lambda_j>0$,  ${\sum}_j\lambda_j=\rho$ and
$$
g_j^TS_\ell g_j\leq {r\kappa\over\varkappa}s_\ell,\;\ell\leq L.
$$
Finally, setting $t=2\ln(4M^2L)$ we have $\varkappa=t\kappa$, and the probability of the event $\varrho(g_j)\leq 1$, $j\leq M$, is at least 1/2.
Thus, we may generate independent $M$-element tuples $g_1,...,g_M$ until a tuple with $\varrho(g_j)\leq1$, $j\leq M$, is generated, and the probability to obtain the desired $\lambda_j's$ and $g_j's$ in $k$ trials is at least $1-2^{-k}$.\qed

\subsection{Proof of Corollary \ref{cor01}}
 The fact that $\bH_\delta$ is a closed convex cone is evident. Under the premise of {\em (i),} for all $j$ there exists $g_j$ such that $h_j=S_\delta g_j$ with $\varrho(g_j)\leq 1$, Consequently,
$$
\Theta={\sum}_j\lambda_jh_jh_j^T=S_\delta^T\underbrace{[{\sum}_j\lambda_jg_jg_j^T]}_{\Xi}S_\delta,
$$
and $(\Xi,\varkappa{\sum}_j\lambda_j)\in\bM_\delta$ by Proposition \ref{prop2}.{\em i,} implying that $(\Theta,\varkappa{\sum}_j\lambda_j)\in\bH_\delta$, as claimed in {\em (i).}
\par
Now, under the premise of {\em (ii)} we have $\Theta=S_\delta^T\Xi S_\delta$ with $(\Xi,\rho)\in \bM_\delta$. By Proposition \ref{prop2}.{\em ii} we can find by a computationally efficient randomized algorithm $\lambda_j\geq0$ and $g_j$, $j\leq M$, such that $\varrho(g_j)\leq1$, $j\leq M$, ${\sum}_j\lambda_j\leq\rho$, such that $\Xi={\sum}_j\lambda_jg_jg_j^T$, implying that $\Theta={\sum}_j\lambda_jh_jh_j^T$ with $h_j=S_\delta g_j$, so that $\pi_\delta(h_j)\leq1$ due to $\varrho(g_j)\leq1$, as claimed in {\em (ii).} \qed

\subsection{Proof of Lemma \ref{lem1}}
Let $h$ satisfy the premise of the lemma. Let us fix $x\in X$ and set $b_i=a_i-Ax$, so that
\[
|h^Tb_i|\leq\varkappa,\;1\leq i\leq n,\;\; \sum_ix_ih^Tb_i=0.
\]
Observe that distribution of $\xi^t_x$ is a ``mixture'' of distributions from $\SG(a_i,\Theta_i)$, $ i\leq n$, i.e., the nature draws at random index $i$ according to the probability distribution $x$ on $\{1,...,n\}$, then draws $\xi^t_i$ from a sub-Gaussian distribution $P^t_i$ which becomes the realization of $\xi^x_t$. Note that  $\xi^t_i$ is sub-Gaussian with parameters $b_i$ and $\Theta_i$ with $\sum_ix_ib_i=0$.
\par
Let $ s >0$, we have
\begin{align*}
\bE\left\{\e^{ s  h^T\xi^t_x}\right\}&=\sum_ix_i\bE\left\{\exp\{ s  h^T\xi^t_i\}\right\}
 \leq \sum_ix_i\exp\left\{s  h^Tb_i+\tfrac{1}{2} s ^2h^T\Theta_i h\right\}\\&=\e^{\tfrac{1}{2} s ^2\sigma^2}\sum_ix_i \e^{s  h^Tb_i}.
\end{align*}
Setting $\beta_i=sh^Tb_i$, we have $\sum_ix_i\beta_i=0$ and $|\beta_i|\leq \beta:=s\varkappa$. So, by convexity of ${\rm e}^x$,
\[
\sum_ix_i\e^{\beta_i}\leq \sum_i x_i\left[\tfrac{\beta+\beta_i}{2\beta}e^{-\beta}+{\tfrac{\beta-\beta_i}{2\beta}e^\beta}\right]=\tfrac{1}{2}\left[e^{-\beta}+e^\beta\right]\leq e^{\beta^2/2},
\]
whence $\sum_ix_i \e^{s  h^Tb_i}\leq\e^{s^2\varkappa^2/2}$ and
$\bE\left\{\e^{ s  h^T\xi^t_x}\right\}\leq \exp\left\{\tfrac{1}{2}s^2(\varkappa^2+\sigma^2)\right\}.
$
Thus, for $\tau\geq 0$ and $s\geq 0$
\begin{align*}
\Prob\{h^T\xi_x\geq\tau\}&\leq \min_{s\geq 0}\bE\left\{\e^{ s N h^T\xi^t_x}\right\}\e^{-Ns\tau}\leq \min_{s\geq 0}\exp\left\{\tfrac{1}{2}Ns^2(\varkappa^2+\sigma^2)-Ns\tau\right\}\\
&\leq \exp\left\{-\frac{\tau^2N}{2(\varkappa^2+\sigma^2)}\right\}.
\end{align*}
Because the same bound holds for  $\Prob\{h^T\xi_x\leq-\tau\}$ (it suffices to replace $h$ with $-h$) this implies \rf{prob000}.\qed
\par
Now, let $g\in\bR^d$ satisfy the premise in (\ref{ifthen}). Then \begin{align*}
\max_{i,j}|g^T[a_i-a_j]|&\leq \varkappa:=\beta/2,\\
\max_ig^T\Theta_ig&\leq \sigma^2:=\beta^2/4,
\end{align*} implying due to \rf{prob000} that
$$
\Prob\{|g^T\xi_x|>1\}\leq 2\exp\{-N/\beta^2\}=\delta.\eqno{\hbox{\qed}}
$$

\end{document}